\documentclass[twoside,a4paper,11pt]{article}
\usepackage{amsmath}
\usepackage{amsthm}
\usepackage{amssymb}
\usepackage{amsfonts}
\usepackage{mathtools}
\usepackage{graphicx}
\usepackage{cite}
\usepackage{fancyhdr}
\usepackage{color}
\definecolor{Blue}{rgb}{0,0,1}
\definecolor{Red}{rgb}{1,0,0}
\usepackage{subfig}
\usepackage[colorlinks=true,linkcolor=Red,urlcolor=Blue,citecolor=Red]{hyperref}
\usepackage[top=3cm, bottom=3cm, left=3cm, right=3cm]{geometry}
\usepackage{tikz}
\usepackage{perpage}
\MakePerPage{footnote}

\theoremstyle{plain}
\newtheorem{theorem}{Theorem}[section]
\newtheorem{proposition}[theorem]{Proposition}
\newtheorem{lemma}[theorem]{Lemma}
\newtheorem{corollary}[theorem]{Corollary}
\theoremstyle{definition}
\newtheorem{definition}[theorem]{Definition}
\newtheorem{example}[theorem]{Example}

\newcommand{\abskeywords}{\par\vspace{0.5cm}\noindent{\bfseries Keywords:}\,\,}
\newcommand{\amssubject}{\par\vspace{0.1cm}\noindent{\bfseries Mathematics Subject Classification [2010]:}\,\,}
\newcommand{\myrule}{\vskip 0.1cm\hrule height 1pt\vskip 0.1cm}
\newenvironment{alphafootnotes}
{\par\edef\savedfootnotenumber{\number\value{footnote}}

\setcounter{footnote}{0}}
{\par\setcounter{footnote}{\savedfootnotenumber}}
\allowdisplaybreaks
\begin{document}

\vspace*{0.05cm}
\begin{center}
{\Large\bfseries Some applications of linear algebraic methods in combinatorics}

 \end{center}

\begin{center}
\begin{alphafootnotes}
{Maryam Khosravi$^{1,}$\footnote{Speaker. Email address: khosravi$_-$m@uk.ac.ir } and Ebadollah S. Mahmoodian$^{2}$}\\[2mm]
{\small $^{1}$Department of Mathematics and Computer Science, Shahid Bahonar University, Kerman, Iran}\\[2mm]
{\small $^{2}$Department of Mathematics, Sharif
University of Technology, Tehran, Iran}
\end{alphafootnotes}
\end{center}
\index{Khosravi, Maryam}
\myrule
\begin{abstract}
In this note, we intend to produce all latin squares from one of them using suitable move which is defined by small trades and do the similar work on 4-cycle systems. These problems, reformulate as  finding basis for the kernel of special matrices, representef to some graphs.
\abskeywords{latin trades, 4-cycle systems, linear independence}
\amssubject{05B30, 15A03.}
\end{abstract}
\myrule

\section{Introduction}\label{sec1}
Linear algebra provides an important tool in combinatorics and graph theory and in the recent decades several articles and books are written about solving some problems in linear algebra and graph theory by using linear algebraic methods and sometimes very serious and surprising results can be drived from simple idea in linear algebra.

An important question, is that can we produce all elements of a class some of combinatorial objects by starting from one of them and defining a special ``move''?

In  this note, we investigate on two different classes ``latin squares" and ``4-cycle systems". Here, we present our moves by trades.

 \section{Main results}\label{sec2}

 \subsection{Latin square}
We start with latine trades. Readers can find some details in \cite{mahdian}.

 \begin{definition}
A \textit{latin square} $L$ of order $n$ is an $n\times n$ array with entries chosen from an $n$-set $N=\{0,1,\dots,n-1\}$, in such a way that each element of $N$ occurs precisely once in each row and in each column of the array. For the ease of exposition, a latin square $L$ will be represented by a set of ordered triples $$\{(i,j;L_{ij}):\ \text{ element } L_{ij}\text{ accures in cell }(i,j)\text{ of the array}\}.$$

 A \textit{partial latin square} $P$ of order $n$ is an $n\times n$ array
in which some of the entries are filled  with the elements from $N$ in such a way that each element of $N$ occurs at most once in each row and at
most once in each column of the array.

For a partial latin square $P$, the shape of P is $S_P=\{(i,j):\ (i,j;P_{ij}\in P\}$, the volume of $P$ is \ $|S_p|$ and
$R^r_P(C^r_P)$ stands for the set of entries occurring in row (column) $r$ of $P$.
\end{definition}

\begin{definition}
A \textit{latin trade} $T=(P,Q)$ of volume $s$ is an ordered set of two partial latin
squares, of order $n$, such that
\begin{enumerate}
\item $S_p=S_Q$.
\item for each $(i,j)\in S_p$, $P_{ij}\neq Q_{ij}$.
\item for each $r$, $0\leq r\leq n-1$, $R_P^r=R^r_Q$ and $C^r_P=C^r_Q$.
\end{enumerate}
\end{definition}
We usually write a latin trade in one table as follows:
$$\begin{array}{|c|c|c|c|c|}
\hline
.&.&2&3&1\\\hline
.&2&.&1&4\\\hline
1&.&0&4&3\\\hline
0&4&1&.&2\\\hline
4&1&3&2&0\\\hline
\end{array}\
\begin{array}{|c|c|c|c|c|}
\hline
.&.&1&2&3\\\hline
.&1&.&4&2\\\hline
4&.&3&1&0\\\hline
1&2&0&.&4\\\hline
0&4&2&3&1\\\hline
\end{array}\quad \Rightarrow\quad
\begin{array}{|c|c|c|c|c|}
\hline
.&.&2_1&3_2&1_3\\\hline
.&2_1&.&1_4&4_2\\\hline
1_4&.&0_3&4_1&3_0\\\hline
0_1&4_2&1_0&.&2_4\\\hline
4_0&1_4&3_2&2_3&0_1\\\hline
\end{array}
$$
A latin trade of volume 4 which is unique (up to isomorphism), is called an
\textit{intercalate.}
Khanban, Mahdian and Mahmoodian \cite{mahdian}, show that for every latin trade $T=(P,Q)$, by applying some suitable intercalates on $P$  we can reach to $Q$. (We shows this actions by addition and this means that every trade can be written as a sum of intercalates.) In this section, we explain the method.

Consider inclusion matrix $M$. Columns of  this matrix correspond to the  $n^3$ ordered triples, of elements  (in lexicographic order) and its rows correspond to the elements of $\{(u_1,u_2)_{\{n_1,n_2\}}:\ 0\leq u_i\leq n-1, \ n_i=1,2,3\}$.
$$M_{(u_1,u_2)_{\{n_1,n_2\}},(x_1,x_2,x_3)}=1\quad\text{ if and only if }u_i=x_{n_i}.$$
Let ${\bf T}$ be a (signed) frequency vector derived from the trade $T=(P,Q)$, i.e.
$${\bf T}_{i,j,k}=\left\{\begin{array}{ll}1\quad&(i,j;k)\in P\\-1&(i,j;k)\in Q\\ 0&\text{otherwise}\end{array}\right.$$
Then $M{\bf T}={\bf 0}$.

\begin{figure}
\centering
\includegraphics[scale=0.2]{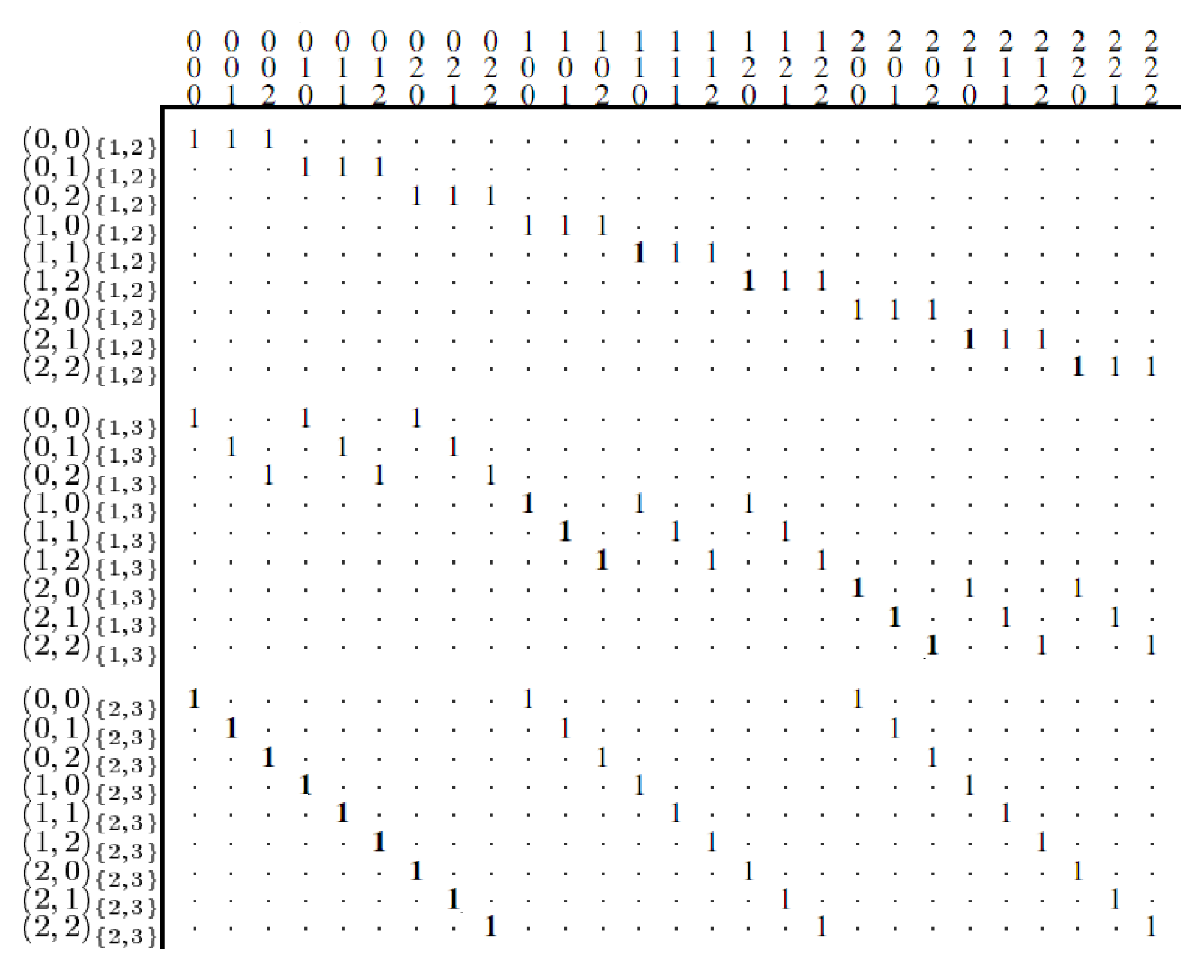}
\caption{Inclusion matrix for $n=3$ }
\end{figure}

 \begin{theorem} If $M$ is the inclusion matrix  with $n^3$ columns, then the nullity of $M$ is $(n-1)^3.$
\end{theorem}
\begin{theorem}
 For all $i,j,k$ the frequency vector of the following intercalates are linearly independent.

\centering {\includegraphics[width=3cm]{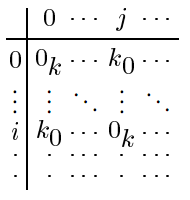}}\end{theorem}
Note that the number of these intercalates is $(n-1)^3$.
Thus we have the following theorem.
\begin{theorem}
Every latin trade can be written as a sum of these intercalates.
\end{theorem}

 \begin{example}
$$\begin{array}{|c|c|c|c|}
\hline
0_1&1_2&2_3&3_0\\\hline
1_2&2_1&.&.\\\hline
2_3&.&3_2&.\\\hline
3_0&.&.&0_3\\\hline
\end{array}=\begin{array}{|c|c|c|c|}
\hline
0_1&1_0&.&.\\\hline
1_0&0_1&.&.\\\hline
.&.&.&.\\\hline
.&.&.&.\\\hline
\end{array}-\begin{array}{|c|c|c|c|}
\hline
0_2&2_0&.&.\\\hline
2_0&0_2&.&.\\\hline
.&.&.&.\\\hline
.&.&.&.\\\hline
\end{array}+\begin{array}{|c|c|c|c|}
\hline
0_2&.&2_0&.\\\hline
.&.&.&.\\\hline
2_0&.&0_2&.\\\hline
.&.&.&.\\\hline
\end{array}$$
$$\qquad-
\begin{array}{|c|c|c|c|}
\hline
0_3&.&3_0&.\\\hline
.&.&.&.\\\hline
3_0&.&0_3&.\\\hline
.&.&.&.\\\hline\end{array}+
\begin{array}{|c|c|c|c|}
\hline
0_3&.&.&3_0\\\hline
.&.&.&.\\\hline
.&.&.&.\\\hline
3_0&.&.&0_3\\\hline
\end{array}
$$
\end{example}

The second combinatoral object that we study, is 4-cycle systems. Some of these result are presented in \cite{khosravi}.
\subsection{4-cycle systems}
Let $G$ be a graph with vertex set $V(G)$ and edge set $E(G)$.
The complete graph of order $n$ is a graph with $n$ vertices for
every pair of them there exists exactly one edge between them and is denoted by $K_n$.

A cycle $(v_1,\ldots,v_n)$ in a graph $G$, is a subgraph of $G$
consists of a set of different vertices $v_i$, and edge set
$\{\{v_1,v_2\},\{v_2,v_3\},\ldots,\{v_n,v_1\}\}$.
\begin{definition}
A 4-cycle system of order $n$, denoted by $4CS(n)$, is a collection  of cycles
of length 4 whose edges partition the edges of $K_n$.
\end{definition}
\begin{definition}
Let $T$ be a subset of edge-disjoint $4$--cycles on vertex set
$\{1,2,\ldots,v\}$. Then $T$ is called a {\sf $4$--cycle trade},
if there exists a set of edge-disjoint $4$--cycles on vertex set
$\{1,2,\ldots,v\}$, $T^{*}$, such that $T^{*}\cap T=\emptyset$ and
$\bigcup_{C\in T}E(C)=\bigcup_{C\in T^{*}}E(C)$.
\end{definition}
 We call $T^{*}$ a {\sf disjoint mate} of $T$ and the pair $(T,T^{*})$ is called a
{\sf $4$--cycle bitrade} of {\sf volume} $t=|T|$ and {\sf
foundation} $v=|\bigcup_{C\in T}V(C)|$.

The following well-known theorem states that for which $n$ a $4CS(n)$ exists.
\begin{theorem}{\cite[Page 266]{8k}}
A necessary
and sufficient condition for the existence of a 4CS(n) is that $n\equiv 1\ (\text{mod }
8)$
\end{theorem}
There are
just four possible configurations of two 4-cycles:
\begin{center}
\definecolor{qqqqff}{rgb}{0.,0.,1.}
\begin{tikzpicture}
\clip(-4.3,-2.48) rectangle (7.06,6.3); \draw (-3.46,5.36)--
(-2.26,5.4); \draw (-2.26,5.4)-- (-2.24,4.34); \draw
(-2.24,4.34)-- (-3.46,4.32); \draw (-3.46,4.32)-- (-3.46,5.36);
\draw (-1.38,5.34)-- (-0.32,5.34); \draw (-0.32,5.34)--
(-0.3,4.38); \draw (-0.3,4.38)-- (-1.4,4.36); \draw (-1.4,4.36)--
(-1.38,5.34); \draw (2.72,5.6)-- (3.5,4.84); \draw (3.5,4.84)--
(2.72,4.); \draw (2.72,4.)-- (1.96,4.76); \draw (1.96,4.76)--
(2.72,5.6); \draw (3.5,4.84)-- (4.3,5.72); \draw (4.3,5.72)--
(5.16,4.96); \draw (5.16,4.96)-- (4.4,4.06); \draw (4.4,4.06)--
(3.5,4.84); \draw (-3.08,2.42)-- (-3.06,1.3); \draw (-3.06,1.3)--
(-1.88,1.3); \draw (-1.88,1.3)-- (-1.88,2.48); \draw
(-1.88,2.48)-- (-3.08,2.42); \draw (-1.88,2.48)-- (-3.06,1.3);
\draw (-3.06,1.3)-- (-1.8,0.34); \draw (-1.88,2.48)--
(-0.74,1.46); \draw (-0.74,1.46)-- (-1.8,0.34); \draw
(3.58,2.42)-- (5.38,1.48); \draw (5.38,1.48)-- (3.64,0.54); \draw
(3.58,2.42)-- (1.96,1.4); \draw (1.96,1.4)-- (3.64,0.54); \draw
(3.64,0.54)-- (3.1,1.48); \draw (3.1,1.48)-- (3.58,2.42); \draw
(3.58,2.42)-- (4.1,1.5); \draw (4.1,1.5)-- (3.64,0.54);
\begin{scriptsize}
\draw [fill=qqqqff] (-3.46,5.36) circle (2.5pt); \draw
[fill=qqqqff] (-2.26,5.4) circle (2.5pt); \draw [fill=qqqqff]
(-3.46,4.32) circle (2.5pt); \draw [fill=qqqqff] (-2.24,4.34)
circle (2.5pt); \draw [fill=qqqqff] (-1.38,5.34) circle (2.5pt);
\draw [fill=qqqqff] (-0.32,5.34) circle (2.5pt); \draw
[fill=qqqqff] (-0.3,4.38) circle (2.5pt); \draw [fill=qqqqff]
(-1.4,4.36) circle (2.5pt); \draw [fill=qqqqff] (2.72,5.6) circle
(2.5pt); \draw [fill=qqqqff] (3.5,4.84) circle (2.5pt); \draw
[fill=qqqqff] (2.72,4.) circle (2.5pt); \draw [fill=qqqqff]
(1.96,4.76) circle (2.5pt); \draw [fill=qqqqff] (4.3,5.72) circle
(2.5pt); \draw [fill=qqqqff] (5.16,4.96) circle (2.5pt); \draw
[fill=qqqqff] (4.4,4.06) circle (2.5pt); \draw [fill=qqqqff]
(-3.08,2.42) circle (2.5pt); \draw [fill=qqqqff] (-3.06,1.3)
circle (2.5pt); \draw [fill=qqqqff] (-1.88,1.3) circle (2.5pt);
\draw [fill=qqqqff] (-1.88,2.48) circle (2.5pt); \draw
[fill=qqqqff] (-1.8,0.34) circle (2.5pt); \draw [fill=qqqqff]
(-0.74,1.46) circle (2.5pt); \draw [fill=qqqqff] (3.58,2.42)
circle (2.5pt); \draw [fill=qqqqff] (5.38,1.48) circle (2.5pt);
\draw [fill=qqqqff] (3.64,0.54) circle (2.5pt); \draw
[fill=qqqqff] (1.96,1.4) circle (2.5pt); \draw [fill=qqqqff]
(3.1,1.48) circle (2.5pt); \draw [fill=qqqqff] (4.1,1.5) circle
(2.5pt);
\end{scriptsize}
\end{tikzpicture}
\end{center}\vspace{-2cm}
It is easily seen that only the last one  forms a trade. We call
this trade a double-diamond.

In \cite{switch,config}, by two different  approaches, the following
theorem was proved.
\begin{theorem}\label{diamond}
For each $n=8k+1$, there exists a $4CS(n)$  which does not have
any double-diamond.
\end{theorem}
The main question in this study is  that can we generate each
$4CS(n)$ from another by using a sequence of double-diamonds? By
Theorem \ref{diamond}, if we restrict the moves to the class of
$4CS(n)$, some of the 4-cycle systems do not generate. Thus, we
need to define the moves in larger space.

To answer this question, we use a linear algebraic approach  as follows.

Let $M$ be a matrix whose rows and columns are corresponded to
edges of $K_n$ and 4-cycles in $K_n$, respectively. Also, for a
trade $(T,T^*)$,  consider the vector $X$ with 1 for cycles in $T$
and $-1$ for cycles in $T^*$. It is easy to see that every vector
$X$ corresponding to a trade is a vector in the kernel of $M$. We
show that the vectors corresponding to double-diamonds form a
basis for $M$.

\begin{lemma}
The matrix $M$ is a full rank matrix.
\end{lemma}
\begin{corollary}
The dimension of kernel of $M$ is $3{n\choose 4}-{n\choose 2}$.
\end{corollary}
\begin{proposition}
For each $n$, there exists $3{n\choose 4}-{n\choose 2}$  linearly
independent vectors in the kernel of $M$, each of them is
corresponding to a double-diamond or a linear combination of
double-diamonds.
\end{proposition}
\begin{proof}
Corresponding to each 4-cycle (except ${n\choose 2}$ of them), we introduce a trade,
containing this 4-cycle, with the following conditions:
\begin{itemize}
\item[1.] The selected 4-cycle does not appear in any of the next trades.
\item[2.] The introduced trade is a double-diamond or a combination
of double-diamonds.
\end{itemize}
Thus, these trades form a triangular basis for the matrix M and
therefore double-diamonds form a basis for the matrix M.
\end{proof}
\begin{corollary}
The vectors corresponding with double-diamonds form a  basis for
$M$. That is, we can construct a 4-cycle system for another one by
using  double-diamonds.
\end{corollary}
Note that in this construction, the way from a 4-cycle  to another
is in a graph $\lambda K_n$ for some $\lambda\in\mathbb{N}$. In
other word, sometimes we have some repeated edge.


 \section{Conclusion}\label{sec4}
Linear algebra is a very strong tool to solve some combinatorial problem and the authors hope that by similar methods, one can find a suitable move between other combinatorial objects.

\end{document}